\documentclass[12pt]{article}
\usepackage{amsmath,amsthm,amssymb,amsfonts}

\newfont{\lie}{eufm10 at 11pt}
\newfont{\liepequenos}{eufm10 at 10pt}

\newfont{\corpos}{msbm10 at 11pt}
\newfont{\corpospequenos}{msbm10 at 10pt}

\newcommand{\Cpequeno}{\mbox{\corpospequenos \symbol{67}}}   
\newcommand{\Hamilpequeno}{\mbox{\corpospequenos \symbol{72}}}
\newcommand{\Rpequeno}{\mbox{\corpospequenos \symbol{82}}}    
\newcommand{\Zpequeno}{\mbox{\corpospequenos \symbol{90}}}    

\newcommand{\C}{\mbox{\corpos \symbol{67}}}          
\newcommand{\Hamil}{\mbox{\corpos \symbol{72}}}      
\newcommand{\R}{\mbox{\corpos \symbol{82}}}          
\newcommand{\Z}{\mbox{\corpos \symbol{90}}}          

\newcommand{\g}{\mbox{\lie g}}       %
\newcommand{\gl}{\mbox{\lie gl}}     

\newcommand{\s}{\mbox{\lie s}}
\newcommand{\pp}{\mbox{\lie p}}
\newcommand{\uni}{\mbox{\lie u}}     %

\newcommand{\sol}{\mbox{\lie so}}

\newcommand{\rr}{\rightarrow}

\newcommand{\hnab}{{\cal H}^{D}}             %
\newcommand{\calV}{{\cal V}}
\newcommand{\calJ}{{\cal J}}

\newcommand{\End}[1]{{\mathrm{End}}\,{#1}}

\newcommand{\Id}{{\mathrm{Id}}}
\newcommand{\dx}{{\mathrm{d}}}

\newcommand{\XIS}{\mbox{\lie X}}   

\newcommand{\na}{\nabla}
\newcommand{\db}{\overline{\partial}}

\newcommand{\vol}{{\mathrm{vol}}}

\newtheorem{teo}{Theorem}[section]

\newtheorem{coro}{Corollary}[section]
\newtheorem{prop}{Proposition}[section]

\numberwithin{equation}{section}

\def\cyclic{\mathop{\kern0.9ex{{+}
\kern-2.2ex\raise-.28ex\hbox{\Large\hbox{$\circlearrowright$}}}}\limits}

\title{Hermitian and quaternionic Hermitian structures on tangent bundles}

\author{Rui Albuquerque\footnote{Departamento de Matem\'atica da Universidade de \'Evora and Centro de Investiga\c c\~ao em Matem\'atica e Aplica\c c\~oes (CIMA), Rua Rom\~ao Ramalho, 59, 7000 \'Evora, Portugal.}\\ rpa@dmat.uevora.pt }

\begin{document}

\maketitle




\begin{abstract}

We review the theory of quaternionic K\"ahler and hyperk\"ahler structures. Then we consider the tangent bundle of a Riemannian manifold $M$ with a metric connection $D$ (with torsion) and with its well estabilished canonical complex structure. With an extra almost Hermitian structure on $M$ it is possible to find a quaternionic Hermitian structure on $TM$, which is quaternionic K\"ahler if, and only if, $D$ is flat and torsion free. We also review the symplectic nature of $TM$. Finally a proper $S^3$-bundle of complex structures is introduced, expanding to $TM$ the well known twistor bundle of $M$.

\end{abstract}

\vspace*{4mm}

{\bf Key Words:} torsion, quaternionic, Hermitian, K\"ahler, symplectic.

\vspace*{2mm}

{\bf MSC 2000:} 53C15, 53C26, 53C55.

\vspace*{10mm}

The author acknowledges the support of Funda\c{c}\~{a}o Ci\^{e}ncia e Tecnologia, through the research center CIMA - Uni\-ver\-si\-dade de \'Evora and through the project POCI/\-MAT/\-{}60671/2004.

\markright{\sl\hfill R. Albuquerque \hfill}

\vspace*{10mm}

\section{Introduction}
\label{Introduction}

The subject of quaternionic Hermitian manifolds still conceals many mysteries for the working geometer. This article starts with a recreation of the main definitions regarding quaternionic K\"ahler structures and their almost immediate properties, pertaining holonomy reduction, which are used later in a particular context.

We develop the theory of complex and quaternionic structures on the tangent bundle of a Riemannian manifold $M$ endowed with a metric connection $D$. It is well known by now how to define an orthogonal almost complex structure $I$ on $TM$ departing from such condition, a construction due to P. Dombrowsky. Such structures have also been studied in a more analytic perspective in \cite{Szoke}. Now, if we assume furthermore that the base manifold is almost Hermitian and take any compatible almost Hermitian $D$, then a sourceful of structures arise on the tangent bundle. We may consider new almost complex structures, orthogonal with respect to the naturally induced metric, as the above $I$, and in a way orthogonal to $I$.

Then $TM$ also carries Hermitian and quaternionic Hermitian structures, and this work concentrates in deciding which conditions on the base space $M$ must be satisfied in order to say wether they are integrable or symplectic and, respectively, quaternionic K\"ahler.

Our techniques involve the determination of the Levi-Civita connection of $TM$ in order to describe the possible holonomy reductions. We hope this is important for other developments of the theory. Our results are confluent with some constructions in \cite{AlekCortesDevch} and the study of quaternionic structures through geometry with torsion is indeed interesting, cf. \cite{Iva}.

\section{Quaternionic K\"ahler structures}
\label{qks}

\subsection{Definitions}
\label{t}

By a quaternionic Hermitian module it is understood a real Euclidian vector space of dimension $4n$ together with a free action by isometries of the Lie group $Sp(1)$ of unit quaternions. This action is assumed to be on the right, as such is the canonical case of $\Hamil^n$. On the Euclidian vector space we also have the left action of $SO(4n)$, which hence contains a copy of the unit quaternions. The automorphisms of the quaternionic Hermitian module constitute another subgroup $Sp(n)\subset SO(4n)$. An isometry $g\in Sp(n)$ if, and only if, $g(vw)=g(v)w$ for any vector $v$ and any $w\in \Hamil$. Hence there is a third resulting subgroup which is the product $Sp(n)Sp(1)$ and which we denote by $G(n)$. Since it is known that the fundamental group of $G(n)$ is $\Z_2$, while $Sp(n)$ is simplyconnected (\cite{Hel}), we have $G(n)=Sp(n)\times_{\Zpequeno_2}Sp(1)$ due to the diagonal action of $\{\pm\Id\}$.

An oriented Riemannian $4n$-manifold $M$ is said to be a \textit{quaternionic K\"ahler} if its holonomy is inside $G(n)$, with an exception in the case $n=1$ -- cf. section \ref{sec:dim4}. If such is the case, then there is a smooth quaternionic Hermitian structure on $M$, i.e. each tangent space $T_xM$ admits a quaternionic Hermitian module structure smoothly varying with $x\in M$. The same is to say $M$ admits a $G(n)$-structure.

Let us reflect upon the implications of the above condition. If the manifold has a $G(n)$-structure this means its frame bundle reduces to a principal $G(n)$-bundle, say $P$. Locally there exist quaternionic Hermitian frames\footnote{These are vector sets $\{v_1,\ldots,v_n\}$ which generate $T_xM$ under right multiplication by scalars in $\Hamilpequeno$. Their existence is proved by the methods in the appendix.} and thus there exists a local lift to an $Sp(n)\times Sp(1)$-structure $\tilde{P}$. The real simple Lie group $Sp(n)$ is the same as $U(2n)\cap Sp(2n,\C)$ (analyze the Lie algebras or simply cf. \cite{Hel}) and hence it has an irreducible representation in $\C^{2n}$, giving rise, locally, to two Hermitian vector bundles:
\begin{equation}\label{EHformalism}
E=\tilde{P}\times_{Sp(n)\times Sp(1)}\C^{2n}\ \ \ \ \ \mbox{and}\ \ \ \ \ H=\tilde{P}\times_{Sp(n)\times Sp(1)}\C^2 
\end{equation}
defined on every sufficiently small open subsets in $M$. One notes $TM\otimes_{\Rpequeno}\C=E\otimes_{\Cpequeno} H$, associated to $P$, in spite of $E,H$ being not, in general, globally defined. Such is known as the $E,H$-formalism\footnote{Nothing as this happens in the geometry of a single almost complex structure, because $GL(1,\Cpequeno)\subset GL(n,\Cpequeno)$.} (cf. \cite{Sal1}).

Recall the metric and the orthogonal complex structure $i1$ in $\C^2$ induce a symplectic 2-form $\omega_H$. Then each $A\in \s\pp(1)=\s\uni(2)=\sol(3)$ is determined by the symmetric 2-product $\omega_H(A\cdot,\cdot)$. In other words, the unit quaternions have Lie algebra (the purely imaginary part of $\Hamil$) a real subspace of the complex vector space $S^2\C^2$, the symmetric complex bilinear forms of $\C^2$. For instance, the unit quaternions $a_11+a_2I+a_3J+a_4K,\ (a_1,\ldots,a_4)\in S^3$ may be represented by taking
\begin{equation}\label{matPauli}
 I=\left[ \begin{array}{cc} 0&1\\-1&0 \end{array}\right] ,\ \ \ \ \ 
   J=\left[ \begin{array}{cc} i&0\\0&-i \end{array}\right],\ \ \ \ \ 
K=IJ=-JI=\left[ \begin{array}{cc} 0&-i\\-i&0 \end{array}\right].
\end{equation} 
Indeed, $I,J,K\in\s\pp(1)$.

As shown, a quaternionic Hermitian structure on a Riemannian manifold does not depend on the complex structure in which we decompose $\Hamil$, but rather on having a real 3-dimensional vector subbundle of $\End{TM}$ over $M$, usually denoted $Q$, locally spanned by three anti-commuting orthogonal almost complex structures ($Q\otimes_{\Rpequeno}\C=S^2 H$). Reciprocally, this induces a $\s\pp(1)\subset \sol(4n)$ associated smooth vector subbundle; hence, by the exponential map, a $Sp(1)$ action on each $T_xM$ smoothly varying with $x$ and therefore a quaternionic Hermitian structure on $M$. We have proved the known result that a $G(n)$ structure is equivalently given by a $Q$ vector bundle as above.

Now the holonomy condition required for a quaternionic K\"ahler manifold corresponds, following the general theory of connections, to the $G(n)$-structure being parallel. The bundle of endomorphisms associated to $\g=\s\pp(n)\oplus\s\pp(1)$ is closed under Levi-Civita covariant differentiation if, and only if, the same happens with the one associated with $\s\pp(1)$, i.e. the rank 3 real vector bundle $Q$. Indeed, notice $\s\pp(n)$ is the centralizer of $\s\pp(1)$ in $\sol(4n)$ and we have
\[ 0=\nabla[\s\pp(n),\s\pp(1)]=[\nabla\s\pp(n),\s\pp(1)]+[\s\pp(n),\nabla\s\pp(1)] .\]
Thus $\nabla\s\pp(n)\subset\s\pp(n)$ if, and only if, $\nabla\s\pp(1)\subset\s\pp(1)$.
\begin{prop}[cf. \cite{Sal1}]
An oriented Riemannian manifold $M$ is quaternionic K\"ahler if, and only if, there exists a parallel vector subbundle $Q\subset\End TM$ locally spanned by three anti-commuting orthogonal almost complex structures.
\end{prop}
As we may check easily, if $q=(I,J,K)$ denotes a \textit{quaternionic triple}, i.e. a local basis of $Q$ of anti-commuting orthogonal almost complex structures, then
\begin{equation}\label{nablaqt}
\nabla q=q\alpha+L
\end{equation}
with $\alpha\in\Omega^1(\s\pp(1))$ and $L_I,L_J,L_K\in\Omega^1(Q^\perp)$ (this is the orthogonal in $\sol(TM)$). Notice $\alpha$ is just a skew-symmetric matrix of 1-forms. The quaternionic K\"ahler condition is thus expressed by the equation $L=0$.

There is also another interesting invariant: any two quaternionic triples $q,q'$  defined on open subsets $U,U'$, respectively, and defining the same structure $Q$ are related by a matrix function $a_{UU'}:U\cap U'\rr SO(3)$, since any $I,J,K$ are pairwise orthogonal and with norm $\sqrt{4n}$. Then in defining the 2-form $\omega_I(X,Y)=\langle IX,Y\rangle$ and $\omega_J,\omega_K$ analogously, we get a well defined 4-form easily seen not to depend on the choice of $q$
\begin{equation}\label{kraines}
\Omega=\omega_I\wedge\omega_I+\omega_J\wedge\omega_J+\omega_K\wedge\omega_K.
\end{equation}
A straightforward computation yields, in the quaternionic K\"ahler case, $\dx\Omega=0$. In general, we find $\dx\Omega=\sum_i\omega_i\wedge\lambda_i$ with the given frame $q$, where $\omega_1=\omega_I,\ \lambda_1=\cyclic_{X,Y,Z}\langle L_I(X)Y,Z\rangle$, etc.

Finally let us recall a third approach to $G(n)$-structures. It is known that $G(n)$ is the set of isometries of a $4n$-dimensional Euclidian vector space for which a non-degenerate $4$-form $\Omega$ defined by (\ref{kraines}) remains invariant (cf. appendix). By a fundamental theorem of Riemannian geometry, the holonomy reduces to $G(n)$ if, and only if, $\nabla\Omega=0$. And it was proved in \cite{Swann} that, when $n>2$, the equation $\dx\Omega=0$ is also a sufficient condition for $G(n)$-holonomy.

\subsection{Topology}
\label{topology}

There is a topological invariant of a quaternionic Hermitian structure which partly measures the obstruction to having globally defined three orthogonal almost complex structures. First notice we have a cohomology sequence associated to 
\begin{equation}\label{sgp}
\begin{array}{ccccccccc}
   1&\rr&\Z_2&\rr& Sp(n)\times Sp(1)&\rr & G(n)  &\rr&1 \\
    &    &    &    &   \downarrow\mathrm{pr}     &     & \downarrow\mathrm{pr}' &    &  \\
   1&\rr&\Z_2&\rr&    Sp(1)         &\rr & SO(3)&\rr&1
\end{array}
\end{equation}
(there exists a projection $\mathrm{pr}'$).
A quaternionic Hermitian structure $P\in H^1(M,G(n))$, a principal $G(n)$-bundle over $M$, lifts to a global principal $Sp(n)\times Sp(1)$-bundle if, and only, if $\delta(P)$ vanishes. The coboundary homomorphism $\delta:H^1(M,G(n))\rr H^2(M,\Z_2)$ follows from the long exact sequence associated to (\ref{sgp}) as a sequence of sheaves of germs of group-valued smooth functions. Recall the second Stiefel-Whitney class $w_2(Q)$ corresponds with the obstruction on lifting the $SO(3)$-structure of $Q=\mathrm{pr}'(P)$ to an $Sp(1)$-structure. Moreover, any $Q\in H^1(M,Sp(1))$ raises to a structure $P$, as explained above through equations (\ref{EHformalism},\ref{matPauli}). We have thus proved $\delta(P)=w_2(Q)$. It measures the existence of $E$ and $H$ globally. 

The picture may be resumed in the following way. Since any two quaternionic triples $q,q'$ defined on open subsets $U,U'$, respectively, are related by a matrix function $a_{UU'}:U\cap U'\rr SO(3)$, a given family of quaternionic triples on an open covering of $M$ gives a cocycle $Q\in H^1(M,SO(3))$; which arrives from a cocycle in $H^1(M,Sp(1))$ if, and only if, $\delta(P)=0$.

\subsection{Hyperk\"ahler and locally hyperk\"{a}hler}
\label{hyperk}

A given Riemannian holonomy is called \textit{hyperk\"ahler} if it reduces from $SO(4n)$ to $Sp(n)\subset G(n)$. In this case the existence of a covering of $M$ by local quaternionic frames with transition functions in $Sp(n)$ only, is implied from the start (in particular $\delta(P)=0$). From this we may construct a global quaternionic triple $I,J,K$ and we observe that $\s\pp(n)=\uni(2n,I)\cap\uni(2n,J)$ (a straightforward computation). Now the equation for holonomy reduction $\na\s\pp(n)\subset\s\pp(n)$ implies reduction to the unitary Lie algebra or simply $\na I\in\uni(2n,I)$ --- which combined with $I^2=-1$ gives $\na I=0$. The same must hold for $J$. Reciprocally, from $\na I=\na J=0$ we arrive to hyperk\"ahler holonomy.

As it is well known, the condition is equivalent to the metric on $M$ being K\"ahler with respect to each almost complex structure.

Some authors immediately attribute the name hyperk\"ahler to a Riemannian manifold with a global quaternionic triple $q=(I,J,K)$ and such that all $\nabla I=\nabla J=\nabla K=0$ (cf. \cite{Besse}). Of course one of the three equations is superfluous.

The term \textit{locally hyperk\"ahler} is reserved for the case when only the reduced holonomy group is inside $Sp(n)$.

\subsection{In dimension 4}
\label{sec:dim4}

In 4 real dimensions we have $Sp(1)Sp(1)=SO(4)$. Hence a Riemannian structure on an oriented manifold $M$ is the same as a quaternionic Hermitian structure.

Every oriented Riemannian 4-manifold $M$ has a unique parallel qua\-ter\-ni\-o\-nic Her\-mi\-ti\-an structure, since any triple $I,J,K$ is identified to an orthonormal basis of the bundle $\Lambda^2_+$ of self-dual two forms and since $\nabla*=*\nabla$. If we select a vector field $U$ with $\Vert U\Vert=1$, then the quaternionic Hermitian module structure on $TM$, with $X_i=\lambda_iU+A_i,\ A_i\in U^\perp,\ i=1,2$, is well known to be given by
\[ X_1\cdot X_2=(\lambda_1\lambda_2-\langle A_1,A_2\rangle)U+\lambda_1A_2+\lambda_2A_1+A_1\times A_2  \]
where $\langle A_1\times A_2,A_3\rangle=\vol(U,A_1,A_2,A_3)$. Notice $\Vert X\cdot Y\Vert=\Vert X\Vert\Vert Y\Vert$. Then any almost complex structure $I=v\cdot:TM\rr TM$ with $v\in U^\perp,\ \Vert v\Vert=1$ and we easily find $\omega_I=U^b\wedge v^b+*(U^b\wedge v^b)$. This picture has led to the construction in \cite{Alb2} of $G_2$-structures on the 7-manifold which is the unit sphere tangent bundle of $M$.

As it was pointed in \cite{Sal1}, we have a lift of a smooth quaternionic Hermitian structure on $M$ to an $Sp(1)\times Sp(1)$-structure if and only if, $M$ is spin. Hence, in this case, $w_2(Q)=w_2(M)=0$.

In view of the above, we finally recall the exception in the definition of quaternionic K\"ahler 4-manifold: a Riemannian structure which is self dual and has the same curvature properties of any other quaternionic K\"ahler structure, namely it is Einstein.

For hyperk\"ahler manifolds we have further strictness: such a 4-manifold is Ricci flat and has flat $\wedge^{\pm}$ bundles. This is a
consequence of having three parallel self-dual 2-forms and hence $R*=*R$, from which $Ric=0$ follows. If locally there exists one parallel unit vector field $U$, then the hyperk\"ahler manifold is itself flat.

\section{$TM$ and its Levi-Civita connection}
\label{TMaiLCc}

Let $M$ be any Riemannian manifold and $D$ any linear metric connection on $M$.

There exists a canonical \textit{vertical} vector field $\xi$ defined on the manifold $TM$:
\begin{equation}\label{cvf}
 \xi_v=v,\qquad\forall v\in TM,
\end{equation} 
under the identification of $\pi^*TM$ with $\calV=\ker(\dx\pi:TTM\rr\pi^*TM)$, where $\pi:TM\rr M$ is the canonical projection. The connection $D$ induces a splitting $TTM=\hnab\oplus\calV$. Moreover, the tautological section $\xi$ carries all the information to produce the splitting. This has already been thoroughly explained in the context of twistor bundles (cf. \cite{Alb1,Obri}) or of the sphere tangent bundle (cf. \cite{Alb2}), where a similar canonical section $\xi$ was defined.

In sum, it follows from the theory that $X\in\hnab\Leftrightarrow(\pi^*D)_X\xi=0$. Essentially, one proves that $\xi$ varies exactly on vertical directions.

Furthermore, for a given vector field $X\in\Omega(TM)={\XIS}_M$ and vector $v\in T_xM$, the vertical part of $\dx X(v)$ is precisely $D_vX$. The theory gives us a projection map $\pi^*D_{\cdot}\xi$ and thus $(\dx X(v))^v=\pi^*D_{\dx X(v)}\xi=(X^*\pi^*D)_vX^*\xi=D_vX$.

Now, we may endow $TM$ with a Riemaniann structure and an induced metric connection denoted $D^*$. Naturally, the metric is defined via the pull-back metric on $\pi^*TM=\calV$ and the isometry $\dx\pi_|:\hnab\rr\pi^*TM$. The decomposition into horizontals and verticals is orthogonal and the metric connection $D^*$, in fact given by $\pi^*D$, preserves this splitting.

Let $R^*=\pi^*R^{D}=R^{\pi^*D}$ denote the curvature tensor of $D^*$. We have $R^*\xi\in\Omega^2(\calV)$. Notice we use $\,\cdot\,^v,\ \,\cdot\,^h$ to denote the vertical and horizontal parts, respectively, of a $TTM$ valued tensor, but the identity $X^h=\dx\pi(X)$ may appear as well.
\begin{teo}
The Levi-Civita connection $\na$ of $TM$ is given by
\begin{equation}\label{lcTM}
\na_XY=D^*_XY-\dfrac{1}{2}R^*_{X,Y}\xi+A_XY+\tau_XY
\end{equation}
where $A,\tau$ are $\hnab$-valued tensors defined by
\begin{equation}\label{AlcTM}
\langle A_XY,Z^h\rangle=\dfrac{1}{2}\langle R^*_{X^h,Z^h}\xi, Y^v\rangle + \dfrac{1}{2}\langle R^*_{Y^h,Z^h}\xi, X^v\rangle
\end{equation}
and
\begin{equation}\label{taulcTM}
\tau(X,Y,Z)=\langle \tau_XY,Z^h\rangle= \dfrac{1}{2}\bigl(T(Y,X,Z)-T(Z,X,Y)+T(Y,Z,X)\bigr),
\end{equation}
with $T(X,Y,Z)=\langle \pi^*T^D(X,Y),Z\rangle$, for any vector fields $X,Y,Z$ over $TM$.
\end{teo}
\begin{proof}
Let us first see the horizontal part of the torsion:
\begin{eqnarray*}
\dx \pi(T^\na(X,Y))& =& D^*_XY^h+A_XY+\tau_XY-D^*_YX^h-A_YX-\tau_YX-\dx\pi[X,Y] \\ 
&=& \pi^*T^D(X,Y)+\tau_XY-\tau_YX,
\end{eqnarray*}
since this is how the torsion tensor of $M$ lifts to $\pi^*TM$ and since $A$ is symmetric. Now we check the vertical part.
\begin{eqnarray*}
(T^\na(X,Y))^v& = & D^*_XY^v-\dfrac{1}{2}R^*_{X,Y}\xi-D^*_YX^v+ \dfrac{1}{2}R^*_{Y,X}\xi-[X,Y]^v\\
& = &D^*_XD^*_Y\xi-R^*_{X,Y}\xi-D^*_YD^*_X\xi-D^*_{[X,Y]}\xi\ =\ 0.
\end{eqnarray*}
$\na$ is a metric connection if, and only if, the difference with $D^*$ is skew-adjoint. This is an easy straightforward computation: on one hand
\begin{equation*}
\langle (\na-D^*)_XY,Z\rangle = -\dfrac{1}{2}\langle R^*_{X,Y}\xi,Z^v\rangle +\dfrac{1}{2}\langle R^*_{X^h,Z^h}\xi, Y^v\rangle + \dfrac{1}{2}\langle R^*_{Y^h,Z^h}\xi, X^v\rangle+\tau(X,Y,Z).
\end{equation*}
On the other hand,
\begin{equation*}
\langle (\na-D^*)_XZ,Y\rangle = -\dfrac{1}{2}\langle R^*_{X,Z}\xi,Y^v\rangle +\dfrac{1}{2}\langle R^*_{X^h,Y^h}\xi, Z^v\rangle + \dfrac{1}{2}\langle R^*_{Z^h,Y^h}\xi, X^v\rangle+\tau(X,Z,Y).
\end{equation*}
hence the condition is expressed simply by $\tau(X,Y,Z)=-\tau(X,Z,Y)$. This, together with $\pi^*T^D(X,Y)+\tau_XY-\tau_YX$, determines $\tau$ uniquely as the form given by (\ref{taulcTM}).
\end{proof}
We remark that from the formula it is clear that $\hnab$ corresponds with an integrable distribuition if, and only if, the Riemmannian manifold $M$ is flat. Indeed, the vertical part of $[X,Y]=\na_XY-\na_YX$, for any pair of horizontal vector fields, is $R^*_{X,Y}\xi$.

Notice $R^*_{X,Y}\xi$ and $\tau(X,Y,Z)$ are null if one of the directions $X,Y,Z$ is vertical. With $A_XY$ the same happens if both $X,Y$ are vertical or horizontal.

It is important to understand when the tensor $\tau$ vanishes. By a result of \'E. Cartan, cf. \cite{Agri}, it is known that the space of torsion tensors $\wedge^2TM\otimes TM$ of a metric connection decomposes into irreducible subspaces like
\begin{equation}\label{decompoftorsion}
{\cal A}\oplus\wedge^3TM\oplus TM,
\end{equation}
where $\wedge^3$ is the one for which $\langle T(X,Y),Z\rangle$ is completely skew-symmetric and where $TM$ is the subspace a vectorial type torsions, i.e. for which there exists $V\in{\XIS}_M$ such that $T(X,Y)=\langle V,X\rangle Y-\langle V,Y\rangle X$. $\cal A$ is an invariant subspace orthogonal to those two. We have the following result:
\begin{prop}\label{tauzero}
$\tau=0$ if, and only if, $T^D=0$.
\end{prop}
\begin{proof}
If $\tau=0$, then $T(Y,X,Z)=T(Z,X,Y)+T(Z,Y,X)$; by the symmetries in $X,Y$ this tensor vanishes.
\end{proof}

\subsection{A complex structure on $TM$}
\label{AcsoTM}

Let $\theta\in\End{TTM}$ be the map which sends $\hnab$ isomorphically onto $\calV$, in view of each subspace $T_vTM$ being identified with $T_{\pi(v)}M\oplus T_{\pi(v)}M$. We see $\theta$ as an endomorphism, imposing $\theta\calV=0$. With respect to the metric we defined above on $TM$ the adjoint of $\theta$ verifies $\theta^t(\calV)=\hnab,\ \theta^t(\hnab)=0$. The following map
\begin{equation}\label{csTM}
 I=\theta^t-\theta 
\end{equation}
is a compatible almost complex structure on $TM$. Indeed, $\theta^t\theta=1_{\hnab}\oplus 0,\ \theta\theta^t=0\oplus1_{\calV}$. For any metric connection, in general, we easily deduce $\nabla\theta^t=(\nabla\theta)^t$ and, for any compatible almost complex structure $I$, \begin{equation}\label{eq.der.cov.omega_I}
\nabla_X\omega_I\,(Y,Z)=\langle(\nabla_X I)Y,Z\rangle.
\end{equation}
For the moment we have $D^*\theta=0$ and hence $D^*I=0$.
\begin{teo}\label{integrabilidadedeI}
(i) The following two assertions are equivalent: $(TM,I)$ is a complex manifold; $D$ is torsion free and flat. If any of these occur, then $M$ is a flat Riemannian manifold and $TM$ is K\"ahler flat.\\
(ii) $\omega_I$ is closed if, and only if, $D$ is torsion free.
\end{teo}
\begin{proof}
On any Riemannian manifold a compatible almost complex structure is integrable if, and only if, $\nabla_uv$ is in the $+i$-eigenbundle of $I$ for all $u,v$ in this same eigenbundle (cf. \cite{Sal2}). The sufficiency of this condition is trivial to prove: if $\nabla_uv$ is in the $+i$-eigenbundle, then the same is true for $[u,v]=\na_uv-\na_vu$. The necessity comes from $\langle[u,v],w\rangle=0$ implying $\langle\na_uv,w\rangle$ to be both a skew- and symmetric 3-\textit{tensor}.

Let us prove (i). In our case, $Iu=iu$ is equivalent to $u=u^h+i\theta u^h$, i.e. the $+i$ eigenbundle $T'TM\simeq\hnab\otimes\C$. Indeed $(\theta^t-\theta)u=-\theta u^h+iu^h=iu$ and the dimensions agree. So we may take $u=X+i\theta X,\ v=Y+i\theta Y$, with $X,Y\in\hnab$ real horizontal vector fields. By (\ref{lcTM})
\begin{eqnarray*}
\na_uv&=&D^*_uv-\frac{1}{2}R^*_{u,v}\xi+A_uv+\tau_uv\\
&=& D^*_uv-\frac{1}{2}R^*_{X,Y}\xi+i(A_X\theta Y+A_{\theta X}Y)+\tau_XY.
\end{eqnarray*}
Now the condition resumes to 
\begin{eqnarray*}
i\theta(i(A_X\theta Y+A_{\theta X}Y)+\tau_XY)=-\frac{1}{2}R^*_{X,Y}\xi.
\end{eqnarray*}
The imaginary part of this gives $\tau=0$ or $T^D=0$ by corollary \ref{tauzero}. For the real part, doing the inner product with a vertical vector gives an equation which we may further simplify by the first Bianchi identity ($D$ is torsion free). It yields the vanishing of the curvature tensor $R^D$. Therefore $\na I=D^*I=0$ and the result follows.

Now we prove (ii) (which implies the second part of (i)). Consider a unitary frame on $TM$ $e_1,\ldots,e_{m},\theta e_1,\ldots,\theta e_{m}$ induced from an orthonormal frame on $M$. Let $e_{i+m}=\theta e_i$. By (\ref{eq.der.cov.omega_I}) and $[R^*_{e_i,\cdot}\xi,\theta]=0$, we have
\begin{eqnarray*}
\dx\omega_I &=& \sum_{i=1}^{2m}\na_i\omega_I\wedge e^i\ =\ \frac{1}{2} \sum_{i,j,k=1}^{2m}\langle\nabla_i(\theta^t-\theta)\,e_j,e_k\rangle e^{ijk}\ =\ -\sum\langle(\nabla_i\theta)e_j,e_k\rangle e^{ijk}\\
& & \ \ \ =\ -\sum\langle (A+\tau)_{e_i}\theta e_j-\theta (A+\tau)_{e_i}e_j,e_k\rangle e^{ijk}.
\end{eqnarray*}
Since $A$ is symmetric and $\tau_{ijk}$ vanishes when $i,j$ or $k$ is vertical, we get $\dx\omega_I=$
\begin{eqnarray*}
= -\sum_{i,j,k=1}^{2m}\langle A_{e_i}\theta e_j-\theta\tau_{e_i}e_j,e_k\rangle e^{ijk}\ =\ \sum_{i,j,k=1}^{m}-\frac{1}{2}\langle R^*_{ik}\xi,\theta e_j\rangle e^{ijk}+\tau_{ijk} e^{ijk+m}\ = \qquad\\
= -\sum_{i<j<k}^m(\langle R^*_{ik}\xi,\theta e_j\rangle -\langle R^*_{jk}\xi,\theta e_i\rangle-\langle R^*_{ij}\xi,\theta e_k\rangle)e^{ijk}+ \sum_{i<j}^m\sum_{k=1}^{m}(\tau_{ijk}-\tau_{jik}) e^{ijk+m}.
\end{eqnarray*}
Since the skew-symmetric part in $X,Y$ of $\tau(X,Y,Z)$ is the torsion of $D$, up to a constant, we must have 0 torsion and thence, by the Bianchi identity, the rest of $\dx\omega_I$ vanishes as well.
\end{proof}
We remark the equivalence in part (i) of the theorem is due to P. Dombrowski, cf. \cite{Dom}, seemingly the first to discover and study the structure $I$.

Notice $\omega_I$ over $TM$ looks very much the same as the natural closed symplectic structure on the co-tangent bundle $T^*M$ of any smooth manifold. Up to the metric-induced isomorphism, we have proved these two are the same if, and only if, we consider the Levi-Civita connection of $M$.

\subsubsection{A remark on complex structures on vector bundles}
\label{Arocsovb}

We recall here some details from the theory of holomorphic vector bundles. Let $M$ be a complex manifold and $E\stackrel{\pi}\rr M$ denote a complex vector bundle of rank $k$, so that it has a smooth complex structure $J=i$. Also let $D$ denote a complex connection on $E$, i.e. one for which $J$ is parallel. 

Recall there exists a natural $\db^E$ operator on sections of $E$ when this is holomorphic.

The following well known result is due to Koszul and Malgrange, cf. \cite{KosMal}. A vector bundle $E$ admits a holomorphic structure such that $\db^Ee=D''e:=\mathrm{pr}\circ De$, where $e$ is any section and $\mathrm{pr}$ is the projection onto the $-i$ eigenbundle $T^*M^{(0,1)}\otimes E$, if, and only if, the $(0,2)$ part of the curvature $R$ of $D$ vanishes. Moreover the holomorphic structure is unique with such condition.

The proof is simple: if we write $E=P\times_{GL(k,\Cpequeno)}\C^k$ with $P$ a principal bundle and use a global $\gl(k,\C)$-valued connection 1-form $\alpha$ to describe $D$ and a local chart $z:U\rr\C^n$ of $M$, then the components of $\alpha$ plus the components of $\pi^*\dx z$ are sufficient to generate a subspace of, imposed, $(1,0)$- $GL(k,\C)$-equivariant forms, and therefore a bundle compatible almost complex structure on $P$, and hence on $E$. By Newlander-Niremberg's celebrated theorem, such structure is integrable if, and only if, the subspace generates a $\dx$-closed ideal in the space of differential forms. This is equivalent to the vanishing of $(\dx\alpha)^{(0,2)}=(\rho-\alpha\wedge\alpha)^{(0,2)}=\rho^{(0,2)}$ where $\rho$ is the curvature form.

The uniqueness of the holomorphic structure \textit{with} the condition $\db^E=D''$ follows, since it is known that it is univocally determined by the underlying almost complex structure and the latter is determined by $\pi$ and $\alpha$ globally.

We may draw a further conclusion: the holomorphic structure of $E$ is the same for all $D$ for which $\rho^{(0,2)}=0$ \textit{and} the connection 1-form is type $(1,0)$, $\alpha''=0$.

We remark that the uniqueness of $D$ is sometimes mistakenly inferred in some of the literature, but it is not even the case in a Hermitian setting as the most trivial example will show; consider $M=\C$ and $D$ nontrivial on the tangent bundle with canonical complex structure, $D=\dx+\mu$, with $\mu$ any $i\R$-valued 1-form. Also $R^D=\db\mu-\partial\overline{\mu}$ is a pure imaginary 2-form which may well not vanish.

In the Hermitian case with {\it the} Hermitian connection, unique as Hermitian and type $(1,0)$ connection, we may say $D$ is flat if, and only if, the connection 1-form is holomorphic. This is because the curvature can only be $(1,1)$, by the metric symmetries, and therefore $\rho=\db\alpha$.

Refering the naturally holomorphic tangent bundle of any complex manifold, furnished with a complex linear connection with $R^{(0,2)}=0$, we have a simple criteria to see if $\db^{TM}=D''$, and reciprocally: the torsion of $D$ must be $(2,0)$. Essentially, this is because the torsion form coincides with $\alpha\wedge\dx z$.

\section{Natural complex structures on $TM$ with almost Hermitian $M$}
\label{CsoTMwaHM}

\subsection{The second complex structure, a pair of them}
\label{Tscsap}

Let $(M,\calJ)$ be an almost Hermitian manifold of real dimension $m=2n$. Let $D$ denote a linear Hermitian connection: a metric connection satisfying $D\calJ=0$. In the following we adopt the notation from the last section.

We may define two natural almost complex structures on $TM$, which we denote by $J$ or $J^\pm$: admiting again the decomposition of $TTM$ into $\hnab\oplus\calV$ we write
\begin{equation}\label{Jmaisemenos}
J^\pm=\calJ\oplus\pm\calJ.
\end{equation}
And let, as usual, $T'M$ denote the $+i$-eigenbundle of $\calJ$.
\begin{teo}\label{integrabilidadedeJ}
(i) $J^+$ is integrable if, and only if, $\calJ$ is integrable and the curvature of $D$ verifies $R^D_{u,v}\overline{w}=0,\ \forall u,v,w\in T'M$.\\
(ii) $J^-$ is integrable if, and only if, $\calJ$ is integrable and $R^D_{u,v}w=0,\forall u,v,w\in T'M$.\\
(iii) $(TM,\omega_{J^\pm})$ is symplectic if, and only if, the Hermitian connection $D$ is flat and its torsion verifies
\begin{equation}\label{eqtorTMJsymp}
 T\in[[{\cal A}]]\oplus[[\XIS_M]].
\end{equation}
This meaning\footnote{We write $[[\cal A]]={\cal A}'+{\cal A}''$ for a vector space of tensors on $T'M$ plus the conjugate of ${\cal A}'$.} that: $T$ has no totally skew-symmetric part, according to (\ref{decompoftorsion}), and $T$ is $(3,0)+(0,3)$ with respect to $\calJ$.
\end{teo}
\begin{proof}
Let $u,v,w$ denote vectors in the $+i$-eigenbundle of $J$. The integrability equation is $(1+iJ)\na_uv=0,\ \forall u,v$. Equivalently, since $\calJ,J$ are $D,D^*$ parallel, respectively, we have
\[ \mbox{(a)}\,\ (1\pm i\calJ)R^*_{u,v}\xi=0\qquad \mbox{and}\qquad \mbox{(b)}\,\ (1+ i\calJ)(A_uv+\tau_uv)=0   \] 
according to vertical and horizontal types. So the two curvature conditions in (i) and (ii) correspond to (a). With respect to (b), in particular for $u,v\in{\hnab}\,'$ we must have $\tau_uv\in{\hnab}\,'$. By a straightforward argument as in corollary \ref{tauzero}, this is the same as $\pi^*T^D(u,v)\in{\hnab}\,'$, or $\pi^*[\pi_*u,\pi_*v]\in{\hnab}\,'$ --- corresponding on the base manifold $M$ to the integrability of $\calJ$. For $u,w$ horizontal and $v$ vertical, since the metric on $M$ is a (1,1) tensor, (b) reads equivalently as $\langle A_uv,w\rangle=0$. Which is
\[ \langle R^*_{u,w}\xi,v\rangle=\tfrac{1}{2}\langle (1\mp i\calJ)R^*_{u,w}\xi,v\rangle=0,\]
due to (a). But this is always true since the projection $\tfrac{1}{2}(1\pm i\calJ)v=0$.

Now let us see assertion (iii). We first compute,
\begin{eqnarray*}
\langle(\na_XJ)Y,Z\rangle &=&\langle-\frac{1}{2}[R^*_{X,\cdot}\xi+A_X+\tau_X,J]Y,Z\rangle\\
&=& \langle -\frac{1}{2}R^*_{X^h,\calJ Y^h}\xi\pm\frac{1}{2}\calJ R^*_{X^h,Y^h}\xi,Z^v\rangle+\langle \pm A_{X^h}\calJ Y^v+\\
& & \qquad +A_{X^v}\calJ Y^h-\calJ A_{X}Y-\calJ\tau_{X^h}Y^h+\tau_{X^h}\calJ Y^h,Z^h\rangle.
\end{eqnarray*}
We denote $R_{\alpha\beta\gamma}=\langle R^*_{e_\alpha,e_\beta}\xi,e_\gamma\rangle$, with $\calJ e_\alpha$ represented by $\hat\alpha$, for an orthonormal frame $e_1,\ldots,e_{m},e_{1+m}=\theta e_1,\ldots,e_{m+m}=\theta e_{m}$ induced from an orthonormal frame of $M$. Now using the symmetry of $A$,
\begin{eqnarray*}
\dx\omega_J &=& \sum_{i=1}^{2m}\na_i\omega_J\wedge e^i\ =\ \frac{1}{2}\sum_{i,j,k=1}^{2m}\langle(\na_iJ)e_j,e_k\rangle e^{ijk}\ = \\
&=& \sum_{i,j,k=1}^{m}-\frac{1}{4}R_{i\hat{j}k+m}e^{ijk+m}\mp\frac{1}{4} R_{ij\widehat{k+m}}e^{ijk+m}\pm \frac{1}{4}R_{ik\widehat{j+m}}e^{i,j+m,k}+\\
& &\qquad\qquad +\frac{1}{4}R_{\hat{j}ki+m}e^{i+m,j,k}+ \frac{1}{2}(\tau_{ij\hat{k}}+\tau_{i\hat{j}k})e^{ijk}\\
&=& \sum_{i,j,k=1}^{m}\frac{1}{4} \bigl(-R_{i\hat{j}k+m}\mp R_{ij\widehat{k+m}}\mp R_{ij\widehat{k+m}}-R_{\hat{j}ik+m}\bigr)e^{ijk+m}+ 
\frac{1}{2} (\tau_{ij\hat{k}}-\tau_{ik\hat{j}})e^{ijk}.
\end{eqnarray*}
Since $\tau_{ijk}$ is skew-symmetric in $j,k$, we get
\begin{equation}\label{domegaJ}
\begin{split}
\dx\omega_J \ =\ \sum_{i,j,k}^{m}\mp \frac{1}{2}R_{ij\widehat{k+m}}e^{ijk+m}+ \tau_{ij\hat{k}}e^{ijk} \ = \hspace{2cm}\\
         \ =\ \sum_{i<j}\sum_k\mp R_{ij\widehat{k+m}}e^{ijk+m}+ 2\sum_{i<j<k}(\tau_{ij\hat{k}}+\tau_{jk\hat{i}}+\tau_{ki\hat{j}})e^{ijk}.
\end{split}
\end{equation}
Now we are in position to prove (iii). To have $\dx\omega_J=0$ the flatness of $D$ is evident; the cyclic sum in $i,j,k$ of $\tau_{ij\hat{k}}$ above implies
\[T_{ji\hat{k}}-T_{\hat{k}ij}+T_{j\hat{k}i}+T_{ik\hat{j}}-T_{\hat{i}jk}+T_{k\hat{i}j}+T_{kj\hat{i}}-T_{\hat{j}ki}+T_{i\hat{j}k}=0 . \]
If $i,j,k$ are indices of three vectors in $T'M$, then we simplify this to
\[T_{jik}-T_{kij}+T_{kji}=0\]
which is the totally skew part of $T$ on $\otimes^3T'M$. If $i,j$ represent vectors in $T'M$ and $k:=\overline{k}$ in $T''M$, then we find
\begin{equation*}
\begin{split}
-T_{ji\overline{k}}+T_{\overline{k}ij}-T_{j\overline{k}i}+T_{i\overline{k}j} -T_{ij\overline{k}}+T_{\overline{k}ij}+T_{\overline{k}ji}-T_{j\overline{k}i}+T_{ij\overline{k}}\ = \hspace{1cm}\\
 T_{ij\overline{k}}-T_{i\overline{k}j}-3T_{j\overline{k}i}\ =\ 0.
\end{split}
\end{equation*}
Equivalently $3T_{j\overline{k}i}=T_{ij\overline{k}}-T_{i\overline{k}j}$ for \textit{all} indices $i,j,k$. In repeating the equation, we deduce $9T_{j\overline{k}i}=3T_{ij\overline{k}}-T_{ji\overline{k}}+T_{j\overline{k}i}$ or
$8T_{j\overline{k}i}=4T_{ij\overline{k}}$. Hence $T_{j\overline{k}i}$ is totally skew-sym\-metric and this same equation says it must be 0.

Taking conjugates, since $T$ is real, we see both $T_{\overline{i}jk}$ and $T_{i\overline{j}\overline{k}}=0$. In particular, the whole skew-symmetric part of the torsion must vanish. This proves the result.
\end{proof}
Notice for the case $J^+$ we see in part (i) of the theorem that the integrability depends on $R^D_{\overline{u},\overline{v}}w=0,\ \forall u,v,w\in T'M$ (the conjugate of the written condition), just like Koszul-Malgrange's theorem prescribes when we see $E=T'M$ with complex structure $J=i$,  cf. section \ref{Arocsovb}. Moreover part (i) is stronger than this celebrated theorem since it does not assume integrability on the base space.

Let $\omega_\calJ$ denote the 2-form on $M$. It is easy to deduce the formula 
\[ \dx\omega_\calJ(X,Y,Z)=\omega_\calJ(T(X,Y),Z)+\omega_\calJ(T(Y,Z),X) +\omega_\calJ(T(Z,X),Y) ,\]
therefore with little extra work we may show that $T$ satisfies condition (\ref{eqtorTMJsymp}) if, and only if, $(M,\omega_\calJ)$ is a symplectic manifold.

The condition found for the torsion in part (iii) is quite interesting if we confront with the ``QKT-connections'' studied in \cite{Iva}; surprisingly those are required to have $T\in\wedge^3$ and to be type (1,2)+(2,1) with respect to $\calJ$.

\subsection{The third complex structure on $TM$}
\label{TtcsoTM}

This work would not be complete if we did not consider the following almost complex structure on the tangent bundle of the Riemannian manifold $M$. Consider the same setting as above and define $J$ to be $J^-$. Consider also the complex structure $I$ from section \ref{AcsoTM}. Then $K=IJ=-JI$ is a new $D^*$-parallel almost complex structure, since $J\theta=-\theta J$, and hence we must do an analysis regarding complex and symplectic geometries just as previously.
\begin{teo}\label{integrabilidadedeK}
(i) The following three are equivalent: $K$ is integrable; $D$ is flat and torsion free; $(M,\calJ)$ is a flat K\"ahler manifold.\\
(ii) $(TM,\omega_K)$ is symplectic if, and only if, $D$ is torsion free. The same is to say $(M,\calJ)$ is K\"ahler.
\end{teo}
\begin{proof}
First we describe $u$ in the $+i$-eigenbundle of $K$. In a decomposition $K(u^h+u^v)=iu^h+iu^v$, this translates in $u^v=i\calJ\theta u^h$. Thence we may write, $T'TM=\{u=X+i\calJ\theta X:\ X\in\hnab\}\otimes\C$. Now the integrability of $K$, as above, is given by $(1+iK)\na_uv=0,\ \forall u,v\in T'TM$. According to types this is simply
\[ \mbox{(a)}\,\ (1+iK)R^*_{u,v}\xi=0\qquad \mbox{and}\qquad \mbox{(b)}\,\ (1+ iK)(A_uv+\tau_uv)=0.\]
Taking $u\in T'TM$ and $v=Y+i\calJ\theta Y$ alike, we get from (a) the equation $(1+iK)R^*_{X,Y}\xi=0$ and so $D$ is flat. From (b) the condition $\tau_XY=0$ follows. Now let us compute $\dx\omega_K$. It could be seen by a formula, $\sum_{i,j,k=1}^{2m}\langle\na_i(J\theta e_j),e_k\rangle e^{ijk}$, but we shall follow the usual proceedre. First,
\begin{eqnarray*}
\lefteqn{ \langle(\na_XK)Y,Z\rangle\ =\ \langle[-\frac{1}{2}R^*_{X,\cdot}\xi +A_X+\tau_X,K]Y,Z\rangle} \\
&=& \frac{1}{2}\langle R^*_{X^h,\theta^t\calJ Y^v}\xi,Z^v\rangle -\frac{1}{2}\langle\theta^t\calJ R^*_{X^h,Y^h}\xi,Z^h\rangle\\ 
& &-\langle A_{X^h}\theta\calJ Y^h+A_{X^v}\theta^t\calJ Y^v +\tau_{X^h}\theta^t\calJ Y^v,Z^h\rangle+\langle\theta\calJ(A_XY+\tau_XY),Z^v\rangle \\
&=& \frac{1}{2}\langle R^*_{X^h,\theta^t\calJ Y^v}\xi+\theta\calJ(A+\tau)_XY,Z^v\rangle +\frac{1}{2}\langle R^*_{X^h,Y^h}\xi,\theta\calJ Z^h\rangle \\
& &  -\frac{1}{2}\langle R^*_{X^h,Z^h}\xi,\theta\calJ Y^h\rangle-\frac{1}{2}\langle R^*_{\theta^t\calJ Y^v,Z^h}\xi,X^v\rangle -\tau(X^h,\theta^t\calJ Y^v,Z^h).
\end{eqnarray*}
Now with the notation of theorem \ref{integrabilidadedeJ}, we have
\begin{eqnarray*}
2\dx\omega_K &=& \sum_{i,j,k=1}^{2m}\langle(\na_iK)e_j,e_k\rangle e^{ijk}  \\
&=& \sum_{i,j,k=1}^{m}\frac{1}{2}R_{i\hat{j}k+m}e^{i,j+m,k+m}+\frac{1}{2} R_{ij\widehat{k+m}}e^{ijk}-\frac{1}{2}R_{ik\widehat{j+m}}e^{ijk}\\
& &\qquad\qquad -\frac{1}{2}R_{\hat{j}ki+m}e^{i+m,j+m,k}- \tau_{i\hat{j}k}e^{i,j+m,k}-\tau_{ij\hat{k}}e^{ijk+m}\\
&=& \sum_{i,j,k=1}^{m}\frac{1}{2} \bigl(R_{i\hat{j}k+m}+ R_{\hat{j}ik+m}\bigr)e^{i,j+m,k+m} +\frac{1}{2}R_{ij\widehat{k+m}}(e^{ijk}-e^{ikj}) \\
& & \qquad\quad -\tau_{i\hat{j}k}(e^{i,j+m,k}-e^{ikj+m})\ =\ \sum  R_{ij\widehat{k+m}}e^{ijk}+2\tau_{i\hat{j}k}e^{ikj+m}.
\end{eqnarray*}
Then by simple computation
\begin{equation}\label{domegaK}
\begin{split}
\dx\omega_K \ =\ \sum_{i<j<k}^{m} (R_{ij\widehat{k+m}}+R_{jk\widehat{i+m}}+R_{ki\widehat{j+m}})e^{ijk}+ \sum_{i<k}\sum_j(\tau_{i\hat{j}k}-\tau_{k\hat{j}i})e^{ikj+m}\\
=\ \sum_{i<j<k}^{m} \cyclic_{ijk}\,R_{ij\widehat{k+m}}e^{ijk}+ 2\sum_{i<k}\sum_jT_{ik\hat{j}}e^{ikj+m} .
\end{split}
\end{equation}
The result now follows easily, since the vanishing of $T$ implies Bianchi identity and already we had $\calJ R^*\xi=R^*\calJ\xi$. Finally if $T=0$ then $D$ is the Levi-Civita connection and so $\calJ$ is integrable and henceforth K\"ahler.
\end{proof}
In some sense, the complex structure $I$ plays a preponderant role. Notice (ii) above is also equivalent to (ii) from theorem \ref{integrabilidadedeI}.

\section{Quaternionic K\"ahler structures on $TM$}
\label{QKsoTM}

In sections \ref{AcsoTM}, \ref{Tscsap} and \ref{TtcsoTM} we saw how to define a quaternionic triple $(I,J,K)$ over the tangent bundle of an almost Hermitian base $(M,\calJ)$ of dimension $m=2n$. In order to decide if it corresponds to true $G(n)$ holonomy, at least in the case $n>2$, we must compute $\dx\Omega$ where $\Omega$ is the 4-form defined in (\ref{kraines}). To start with, let 
\[ e_1,\ldots,e_n,e_{n+1},\ldots,e_{2n},e_{2n+1},\ldots,e_{3n},e_{3n+1},\ldots,e_{4n}\]
be a frame on $TM$ induced from a unitary frame of $M$: $e_{l+n}=\calJ e_l$, $e_{2n+i}=\theta e_i$, with $1\leq l\leq n,\ 1\leq i\leq 2n$. Then it is easy to deduce
\[ \omega_I=-\sum e^{i,i+2n},\qquad  \omega_J=\sum e^{l,l+n}-e^{l+2n,l+3n},\qquad \omega_K=e^{l+n,l+2n}-e^{l,l+3n}.  \]
\begin{teo}\label{chato}
$(TM,I,J,K)$ is a quaternionic K\"ahler manifold if, and only if, $D$ is flat and torsion free.
\end{teo}
\begin{proof}
In the proof of theorem \ref{integrabilidadedeI} we computed $\dx\omega_I$. Using this and formulae (\ref{domegaJ}) and (\ref{domegaK}) we deduce
\begin{eqnarray*}
\frac{1}{2}\dx\Omega &=& \dx\omega_I\wedge\omega_I+ \dx\omega_J\wedge\omega_J+
\dx\omega_K\wedge\omega_K\\
&=& \sum_{i<j<k}^{2n}\sum_{l=1}^{n}\cyclic_{ijk} \biggl( R_{ijk+2n}(e^{ijkll+2n}+e^{ijkl+n,l+3n})+  \biggr.\\
& & \biggl. \quad +2\tau_{ij\hat{k}}(e^{ijkll+n}-e^{ijkl+2n,l+3n})+R_{ij\widehat{k+m}}(e^{ijkl+n,l+2n}-e^{ijkll+3n})\biggr)+ \\
& &  \biggl. \ \ +\sum_{i<j}^{2n}\sum_{k=2n+1}^{4n}\sum_{l=1}^{n} \biggl(
2T_{ijk-2n}(e^{ijk,l,l+2n}+e^{ijk,l+n,l+3n})+  \biggr.\\
& & \biggl. +R_{ij\hat{k}}(e^{ijk,l,l+n}-e^{ijk,l+2n,l+3n})
+2T_{ij\widehat{k-2n}}(e^{ijk,l+n,l+2n}-e^{ijk,l,l+3n})\biggr)
\end{eqnarray*}
with notation given previously. It is easy to check $\dx\Omega=0$ implies $R^D=0, T^D=0$.
\end{proof}

\subsection{A family of quaternionic K\"ahler structures on $TM$.}
\label{AfoqksoTM}

Here we assume we have a $4n$ manifold endowed with a quaternionic triple $q=(\calJ_1,\calJ_2,\calJ_3)$; we are going to extend these endomorphisms to $TTM$ in a canonical fashion as it was done in section \ref{Tscsap}, but now with a certain connection $D$ known as the Obata connection. The following seems not to be so well known, hence we give a proof.
\begin{prop}[Obata]\label{parallelqtriple}
For every quaternionic Hermitian structure $\theta=(I,J,K)$ there is a metric connection $D$ such that $DI=DJ=DK=0$.
\end{prop}
\begin{proof}
Let $\nabla$ denote any metric connection and let $A_E=(\nabla E)E$, for any $E\in \End{TM}$. Then we have $[A_J,J]=(\nabla J)J^2-J(\nabla J)J=-\nabla J+(\nabla J)J^2=-2\nabla J$, proving we can always find a Hermitian connection: $(\nabla+\frac{1}{2}A_J)J=0$. It is easy to see that $A_J$ is an $\sol(TM)$-valued 1-form. We also have
\[  [A_J,I] =(\na J)JI-I(\na J)J= -(\na J)K+K\na J = [K,\na J]  \]
and hence, letting $D=\na+\frac{1}{4}(A_I+A_J+A_K)$, we find
\begin{eqnarray*}
 DI&=&(\na +\tfrac{1}{2}A_I)I-\tfrac{1}{4}[A_I,I]+\tfrac{1}{4}[A_J,I]+\tfrac{1}{4}[A_K,I]\\
&=&\tfrac{1}{4}\bigl(2\na I+K\na J-(\na J)K-J\na K+(\na K)J\bigr)\\
&=&\tfrac{1}{4}\bigl(2\na I+\na (KJ)-\na (JK)\bigr)\ =\ 0.
\end{eqnarray*}
The same equation holds for $J$ and $K$.
\end{proof}
Now let $I_0=I$ be the endomorphism defined in \ref{AcsoTM} and let
\begin{equation}\label{novosIs}
 I_i=\calJ_i\oplus-\calJ_i,\qquad\forall i=1,2,3,
\end{equation}
as the case $J^-$ in \ref{Tscsap}. Notice $I_3\neq I_1I_2=-I_2I_1$. However, the whole four $I_i$ anti-commute with each other. Hence, for each point $(a,b)\in V^4_2$, the Stiefel manifold of pairs of orthonormal vectors $a,b\in\R^4$, we have a quaternionic triple $(I_a,I_b,I_{a,b})$ given by
\begin{equation}\label{novosIsgerais}
 I_x=x_0I_0+x_1I_1+x_2I_2+x_3I_3,\quad\forall x=a,b,\quad\mbox{and}\quad I_{a,b}=I_aI_b.
\end{equation}
It is easy to verify $I_x^2=-1$ and $I_aI_b=-I_bI_a$. Also we let $\Omega_{a,b}=\omega_a^2+\omega_b^2+\omega_{a,b}^2$  where $\omega_a(X,Y)=\langle I_aX,Y\rangle$, etc.

We then have two extreme examples: $a=(1,0,0,0),\ b=(0,1,0,0)$ yield the case with which we started this section. Theorem \ref{chato} gives further information about $\Omega$.

With $a=(0,1,0,0),\ b=(0,0,1,0)$ we have the other case, where the requirement of a quaternionic Hermitian base $M$ is unavoidable. We have also done the computations of the respective $\dx\Omega_{a,b}=0$ and the condition found was the same as for the first case: the very strict torsion free and flat metric connection $D$. The proof is very much alike using a quaternionic frame. Finally, due to the fact that every $a\in S^3$ is connected by a curve $e^{itx}e^{jty}e^{ktz}$ in $\Hamil$ to $(1,0,0,0)$, it may be possible to prove that theorem \ref{chato} holds for every $(a,b)\in V^4_2$.

Recall that for every almost quaternionic Hermitian manifold $(M,Q=<q>)$, there is an associated twistor space ${\cal Z}(M)\subset Q$, an $S^2$-bundle of endomorphisms $a\calJ_1+b\calJ_2+c\calJ_3$, with $(a,b,c)\in S^2$, defining complex structures in each $T_xM$. Thus we have obtained a ``Hopf-twistor'' extension of such bundle associated to the tangent bundle.

\subsection{Over a Riemann surface $M$}
\label{OaRsM}

In order to speak of quaternionic K\"ahler structures on the tangent bundles of Riemannian manifolds, the cases $n=1$ and 2 are missing. We concentrate on the case $n=1$ and recall the desired condition now is the metric on $TM$ to be self-dual and Einstein.

Let $\xi$ be the canonical vector field (\ref{cvf}) and let $\eta$ be the unit vertical vector field such that $\{\frac{\xi}{c},\eta\}_u=\{\frac{u}{c},\eta_u\}$ is a direct orthonormal basis of $T_{\pi(u)}M$, $\forall u\in TM$, with $c_u=\Vert\xi_u\Vert=\Vert u\Vert$. Let $D$ be the usual metric connection on $M$ and denote 
by $k$ the function $k(u)=\langle R_{\frac{u}{c},v}\frac{u}{c},v\rangle$. We may also write $k=\frac{1}{c^2}\langle R_{\xi_h,\eta_h}\xi,\eta\rangle$
where $\xi_h,\eta_h$ are such that their images under $\theta$ are $\xi,\eta$, respectively, $\theta$ being the map introduced in \ref{AcsoTM}. Suppose the torsion of $D$ is such that
\[ T(\xi,\eta)=f_1\xi+f_2\eta \]
with $f_1,f_2$ real functions. Then the tensor defined in (\ref{taulcTM}) satisfies
\[  \tau_{\cdot}\xi_h=(f_1\xi_h^b+f_2\eta_h^b)\eta_h,\qquad \tau_{\cdot}\eta_h=-\frac{1}{c^2}(f_1\xi_h^b+f_2\eta_h^b)\xi_h . \]
A straightforward computation yields the following formulae for the Levi-Civita connection of $TM$:
\begin{equation}\label{nablaTMSuperfRiemann}
\begin{array}{llll}
\na_{\xi}\xi=\xi & \na_{\xi}\eta=0 & \na_{\xi}\xi_h=\xi_h & \na_{\xi}\eta_h=0 \\
\na_{\eta}\xi=\eta & \na_{\eta}\eta=-\frac{\xi}{c^2} & \na_{\eta}\xi_h=(1+\frac{k}{2}c^2)\eta_h & \na_{\eta}\eta_h=-(\frac{1}{c^2}+\frac{k}{2})\xi_h \\
\na_{\eta_h}\xi=0 & \na_{\eta_h}\eta =-\frac{k}{2}\xi_h & \na_{\eta_h}\xi_h= \frac{k}{2}c^2\eta+f_2\eta_h & \na_{\eta_h}\eta_h=-\frac{1}{c^2}f_2\xi_h \\
\na_{\xi_h}\xi=0 & \na_{\xi_h}\eta=\frac{k}{2}c^2\eta_h & \na_{\xi_h}\xi_h= f_1c^2\eta_h & \na_{\xi_h}\eta_h=-\frac{k}{2}c^2\eta-f_1\xi_h 
\end{array}
\end{equation}
From these and other identities such as $[\xi_h,\eta_h]=-c^2k\eta-f_1\xi_h-f_2\eta_h$ (the most relevant between the Lie bracket computations) we may compute the Riemannian curvature of $TM$. Notice $k,f_1,f_2$ only depend on $x\in M$. The upshot of these calculations is the following result: $TM$ is Einstein if, and only if, $k=0$ and
\[  c^2\eta_h(f_1)-\xi_h(f_2)-c^2f_1^2-f_2^2=0,  \]
still an intriguing equation. In particular, we may conclude with a corollary when $k$ is the Gauss curvature.
\begin{coro}
For a Riemann surface $M$, $TM$ with its canonical metric is an Einstein manifold if, and only if, the Riemannian curvature of $M$ is 0.
\end{coro}

\section{Appendix}

We prove here that $Sp(n)Sp(1)=G(n)$ is the (isotropy) subgroup of $SO(4n)$ which leaves invariant the 4-form $\Omega$ defined by the identity (\ref{kraines}) on an Euclidian $4n$-vector space $V$.

The group $Sp(n)$ is by definition the subgroup of isometries of $V$ which commute with the given quaternionic triple $q=(I,J,K)$.

If $g\in G(n)$ then $\forall X\in V,\ w\in\Hamil,\ g(Xw)=g(X)w'=g(X)ww''$ for some $w''\in S^3\subset\Hamil$. This is, $g$ preserves the quaternionic lines and reciprocally. Hence, to see $g^*\Omega=\Omega$ we are bound to prove it for $g\in Sp(1)$. Immediately we deduce
\[  I^*\omega_I=\omega_I,\qquad I^*\omega_J=-\omega_J,\qquad I^*\omega_K=-\omega_K. \]
Since $I^*(\omega\wedge\omega)=I^*\omega\wedge I^*\omega$ and since all the same is true for $J,K$, we see 
\[ I^*\Omega=J^*\Omega=K^*\Omega=\Omega. \]
To prove the reciprocal we need a lemma: if $Y,Y_1,Y_2,Y_3$ is an orthonormal set such that $\Omega(Y,Y_1,Y_2,Y_3)=4$, then $Y_j\in\mathrm{span}\{IY,JY,KY\}$, $\forall j=1,2,3$. Proof: let
$Y_j=\alpha_jIY+\beta_jJY+\gamma_jKY+Z_j$ with $Z_j$ orthogonal to the quaternionic line spanned by $Y$. Then it is easy to compute from identity (\ref{kraines})
\[ 4=\Omega(Y,Y_1,Y_2,Y_3)=4\det\left[ \begin{array}{ccc}
\alpha_1 & \beta_1 & \gamma_1 \\
\alpha_2 & \beta_2 & \gamma_2 \\
\alpha_3 & \beta_3 & \gamma_3 
\end{array}\right].  \]
But since the $Y_j$ are orthonormal, $\alpha_j^2+\beta_j^2+\gamma_j^2+|Z_j|^2=1$. Now these two equations yield $Z_j=0$, proving the lemma.

Finally, suppose $g\in SO(4n)$ and $g^*\Omega=\Omega$. Then take any quaternionic line, with an orthonormal basis $X,X_1,X_2,X_3$. We want to see the $Y_i=g(X_i)$ are all in the same quaternionic line. Since $\Omega(Y,Y_1,Y_2,Y_3)=\Omega(X,X_1,X_2,X_3)=4$, the lemma gives the result.

\end{document}